\documentclass[amssymb, 11pt]{amsart}
\usepackage{amsmath}
\usepackage{amsthm}
\usepackage{amsfonts}

\newcommand{\yy}{\mathbb{Y}}

\newcommand{\qb}[2]{{\left [{#1 \atop #2} \right]}}
\newlength{\standardunitlength}
\newtheorem{prop}{Proposition}[section]

\newtheorem{lemma}[prop]{Lemma}
\newtheorem{cor}[prop]{Corollary}
\newtheorem{theorem}[prop]{Theorem}

\begin{document}

\begin{center}
\title [Separation cutoffs for random walk] {Separation cutoffs for
random walk on irreducible representations}
\end{center}

\author{Jason Fulman}
\address{University of Southern California\\ Los
Angeles, CA 90089-2532}
\email{fulman@usc.edu}



\thanks{The author was partially supported by NSF grant DMS-050391.}

\thanks{Version of March 8, 2007}

\begin{abstract} Random walk on the irreducible representations
of the symmetric and general linear groups is studied. A separation
distance cutoff is proved and the exact separation distance
asymptotics are determined. A key tool is a method for writing the
multiplicities in the Kronecker tensor powers of a fixed
representation as a sum of non-negative terms. Connections are made
with the Lagrange-Sylvester interpolation approach to Markov
chains. \end{abstract}

\maketitle

\section{Introduction}

The study of convergence rates of random walk on a finite group is a
rich subject; three excellent surveys are \cite{Al}, \cite{Sal} and
\cite{D1}. In recent papers (\cite{F1},\cite{F3},\cite{F6}), the
author posed and studied a dual question, namely what can be said
about the convergence rate of random walk on $Irr(G)$, the set of
irreducible representations of a finite group $G$.

To define the random walk on $Irr(G)$, let $\eta$ be a (not
necessarily irreducible) representation of $G$ whose character is real
valued. From an irreducible representation $\lambda$, one transitions
to the irreducible representation $\rho$ with probability
\[ K(\lambda,\rho):= \frac{d_{\rho} m_{\rho}(\lambda \otimes \eta)}{ d_{\lambda} d_{\eta}}.\] Here $d_{\lambda}$ denotes the dimension of $\lambda$ and $m_{\rho}(\lambda \otimes \eta)$ denotes the multiplicity of
$\rho$ in the tensor product (also called the Kronecker product) of
$\lambda$ and $\eta$. Whereas random walk on $G$ has the uniform
measure as a stationary distribution, random walk on $Irr(G)$ has the
Plancherel measure as a stationary distribution. The Plancherel
measure assigns a representation $\lambda$ probability
$\frac{d_{\lambda}^2}{|G|}$.

There are many motivations for the study of these random walks:
six motivations (with literature references and discussion) appear in
the introduction of the recent paper \cite{F6}. There is no need to
repeat the discussion here, but let us just mention that similar
processes have been studied for compact Lie groups and Lie algebras
(\cite{ER},\cite{BBO}), and that the decomposition of iterated tensor
products of finite groups has been studied by combinatorialists
(\cite{F3},\cite{GC},\cite{GK}). Moreover, these random walks arise in
quantum computing (\cite{MR},\cite{F6}) and have been used to derive
the first error bounds in limit theorems for the
distribution of random character ratios (\cite{F1}, \cite{F4}). As is
clear from \cite{F6} and this paper, these random walks are also a
tractable testing ground for results in Markov chain theory.

To illustrate a result in this paper, we discuss the case of random
walk on $Irr(S_n)$. To do this, recall that two commonly used
distances between probability distributions $P,Q$ on a finite set $X$
are total variation distance
\[ ||P-Q|| := \frac{1}{2} \sum_{x \in X} |P(x)-Q(x)| \] and
separation distance \[ s(P,Q) := \max_{x \in X} \left[ 1 -
\frac{P(x)}{Q(x)} \right].\] The following recent result gave a sharp
total variation distance convergence rate estimate.

\begin{theorem} \label{TVbound} (\cite{F6}) Let $G$ be the symmetric group
$S_n$ and let $\pi$ be the Plancherel measure of $G$. Let $\eta$ be
the $n$-dimensional defining representation of $S_n$. Let $K^r$ denote
the distribution of random walk on $Irr(G)$ after $r$ steps, started
from the trivial representation.
\begin{enumerate}
\item If $r=\frac{1}{2}n\log(n)+cn$ with $c \geq 1$ then \[
||K^r-\pi|| \leq \frac{e^{-2c}}{2}. \]
\item If $r=\frac{1}{2}n\log(n)-cn$ with $0 \leq c \leq \frac{1}{6}
\log(n)$, then there is a universal constant $a$ (independent of $c,n$)
so that \[ ||K^r-\pi|| \geq 1 -ae^{-4c}.\]
\end{enumerate}
\end{theorem}

This paper gives precise separation distance asymptotics. Letting
$s(r)$ denote the separation distance after $r$ steps, it will be
shown that for the walk in Theorem \ref{TVbound} and $c$ fixed in
$\mathbb{R}$, \[ s(n \log(n)+cn) = 1-e^{-e^{-c}}(1+e^{-c}) + O \left(
\frac{\log(n)}{n} \right). \] This expression goes to $0$ as $c
\rightarrow \infty$ and to $1$ as $c \rightarrow -\infty$ and a cutoff
(defined precisely in Section \ref{prelim}) occurs since $cn=o(n
\log(n))$. Note that whereas the total variation cutoff occurs at time
$\frac{1}{2}n \log(n)$, the separation distance cutoff occurs at time
$n \log(n)$. The proof of the separation distance asymptotics (and
also the corresponding result for $GL(n,q)$) consists of two steps:

\begin{enumerate}
\item One must determine at which element $\lambda$ of $Irr(G)$ the
separation distance is attained. This is equivalent to finding the
$\lambda$ which minimizes
$\frac{m_{\lambda}(\eta^r)}{d_{\lambda}}$. This is tricky since the
usual formula for multiplicities in tensor products involves character
values and so both positive and negative terms. Our solution is to
give a subtle rewriting of this expression as a sum of non-negative
terms.

\item Once one knows at which representation the separation distance
is attained, one needs a formula for the separation distance. For the
cases in this paper we indicate how to do this using combinatorial
arguments and the diagonalization (i.e. eigenvalues and eigenvectors)
of the random walk on $Irr(G)$.
\end{enumerate}

A rather remarkable fact is that for the cases studied in this paper,
one can use Lagrange-Sylvester interpolation to carry out Step 2
knowing only the eigenvalues (and not the eigenvectors) of random walk
on $Irr(G)$. A similar trick had been usefully applied in the
one-dimensional setting of birth-death chains (\cite{Br}, \cite{DSa}),
but the state spaces $Irr(S_n)$ and $Irr(GL(n,q))$ are
high-dimensional so it is interesting that the trick can be extended
to this context. A sequel will treat combinatorial examples where
similar ideas can be applied.

The organization of this paper is as follows. Section \ref{prelim}
gives background from Markov chain theory and recalls the
diagonalization of the Markov chain on $Irr(G)$. Section
\ref{symirrep} derives separation distance asymptotics for random walk
on $Irr(S_n)$ when $\eta$ is the defining representation (whose
character on a permutation is the number of fixed points). Section
\ref{glirrep} obtains separation distance asymptotics for random walk
on $Irr(GL(n,q))$ when $\eta$ is the representation whose character
is $q^{d(g)}$, where $d(g)$ is the dimension of the fixed space of
$g$. Section \ref{lagrange} discusses Lagrange-Sylvester
interpolation, giving eigenvector-free proofs of some results of
Sections \ref{symirrep} and \ref{glirrep}.

\section{Preliminaries} \label{prelim}

This section collects some background on finite Markov chains, using
random walk on $Irr(G)$ as a running example. Let $X$ be a finite set
and $K$ a matrix indexed by $X \times X$ whose rows sum to 1. Let
$\pi$ be a distribution such that $K$ is reversible with respect to
$\pi$; this means that $\pi(x) K(x,y) = \pi(y) K(y,x)$ for all $x,y$
and implies that $\pi$ is a stationary distribution for the Markov
chain corresponding to $K$.

As an example, the Markov chain on $Irr(G)$ defined in the
 introduction is reversible with respect to the Plancherel measure
 $\pi$. To see this let $\chi$ denote the character of a
 representation, and recall the formula \[ m_{\rho}(\lambda \otimes
 \eta) = \frac{1}{|G|} \sum_{g \in G} \chi^{\lambda}(g) \chi^{\eta}(g)
 \overline{\chi^{\rho}(g)}. \] The equation
 $\pi(\lambda)K(\lambda,\rho)=\pi(\rho) K(\rho,\lambda)$ follows
 because $\eta$ was assumed to be real valued; in fact this was the
 reason for imposing this condition on $\eta$.

 Define $\langle f,g \rangle = \sum_{x \in X} f(x) g(x) \pi(x)$ for
real valued functions $f,g$ on $X$, and let $L^2(\pi)$ denote the
space of such functions. Then when $K$ is considered as an operator on
$L^2(\pi)$ by
\[ Kf(x) := \sum_y K(x,y) f(y),\] it is self adjoint. Hence $K$ has an
orthonormal basis of eigenvectors $f_i(x)$ with $Kf_i(x) = \beta_i
f_i(x)$, where both $f_i$ and $\beta_i$ are real. It is easily shown
that the eigenvalues satisfy $-1 \leq \beta_{|X|-1} \leq \cdots \leq
\beta_1 \leq \beta_0=1$. One calls $K$ ergodic if
$|\beta_{|X|-1}|,|\beta_1|<1$.

As an example, Lemma \ref{diaggroup} determines an orthonormal basis
of eigenvectors for the Markov chains on $Irr(G)$.

\begin{lemma} \label{diaggroup} (\cite{F2}, Proposition 2.3) Let $K$ be the Markov chain on $Irr(G)$ defined using a representation $\eta$ whose character is real valued. The eigenvalues of $K$ are indexed by conjugacy classes $C$ of $G$:
\begin{enumerate}
\item The eigenvalue parameterized by $C$ is
$\frac{\chi^{\eta}(C)}{d_{\eta}}$.
\item An orthonormal basis of eigenfunctions $f_C$ in $L^2(\pi)$ is
defined by $f_C(\rho) = \frac{|C|^{1/2} \chi^{\rho}(C)}{d_{\rho}}$.
\end{enumerate}
\end{lemma}

For instance when $G=S_n$ and $\eta$ is the n-dimensional defining
representation, the eigenvalues are $\frac{i}{n}$ where $0 \leq i \leq
n-2$ or $i=n$, with multiplicity equal to the number of conjugacy classes
of permutations with $i$ fixed points. Similarly, suppose that
$G=GL(n,q)$ and that $\eta$ is the representation whose character is
the number of fixed vectors of $g$ in its natural action on the
n-dimensional vector space $V$. Then the eigenvalues are $q^{-i}$ for
$i=0,\cdots,n$, with multiplicity equal to the number of conjugacy
classes of elements of $GL(n,q)$ with an $n-i$ dimensional fixed space.

A common way to quantify convergence rates of Markov chains is using
total variation distance. Given probabilities $P,Q$ on $X$, one
defines the total variation distance between them as
$||P-Q||=\frac{1}{2} \sum_{x \in X} |P(x)-Q(x)|$. It is not hard
to see that \[ ||P-Q|| = \max_{A \subseteq X} |P(A)-Q(A)| .\] Let
$K_x^r$ be the probability measure given by taking $r$ steps from the
starting state $x$. Researchers in Markov chains are interested in the
behavior of $||K_x^r - \pi||$.

Lemma \ref{genbound} relates total variation distance to the spectrum
of $K$. Part 1 is the usual method for computing the power of a
diagonalizable matrix. Part 2 is proved in \cite{DH} and upper bounds
$||K_x^r - \pi||$ in terms of eigenvalues and eigenvectors and is
effective in many examples; it was crucial in the proof of Theorem
\ref{TVbound} stated in the introduction.

\begin{lemma} \label{genbound}
\begin{enumerate}
\item $K^r(x,y) = \sum_{i=0}^{|X|-1} \beta_i^r f_i(x) f_i(y) \pi(y)$ for any $x,y \in X$.
\item  \[ 4 ||K_x^r - \pi||^2 \leq \sum_y \frac{|K^r(x,y) - \pi(y)|^2}{\pi(y)} = \sum_{i=1}^{|X|-1} \beta_i^{2r} |f_i(x)|^2 .\] Note that the final sum does not include $i=0$.
\end{enumerate}
\end{lemma}

Another frequently used method to quantify convergence rates of Markov
chains is to use separation distance, introduced in
\cite{AD1},\cite{AD2}. The separation distance between probabilities
$P,Q$ on $X$ is defined as \[ s(P,Q) = \max_{x \in X} \left[ 1 -
\frac{P(x)}{Q(x)} \right].\] Since $||P-Q|| = \sum_{x: Q(x) \geq P(x)}
[Q(x)-P(x)]$, it is straightforward that $||P-Q|| \leq
s(P,Q)$. Specializing to random walk on $Irr(G)$ started at the
trivial representation $\hat{1}$, one has that \[ s(K_{\hat{1}}^r,\pi)
= \max_{\lambda} \left[ 1 - \frac{|G|
K^r(\hat{1},\lambda)}{d_{\lambda}^2} \right].\] Lemma 3.2 of \cite{F6}
gives that $K^r(\hat{1},\lambda) = \frac{d_{\lambda}}{(d_{\eta})^r}
m_{\lambda}(\eta^r)$. Thus the separation distance is attained at the
$\lambda$ which minimizes $\frac{m_{\lambda}(\eta^r)}{d_{\lambda}}$.

Finally, let us give a precise definition of the cutoff phenomenon. A
nice survey of the subject is \cite{D2}; we use the definition from
\cite{Sal}. Consider a family of finite sets $X_n$, each equipped with
a stationary distribution $\pi_n$, and with another probability
measure $p_n$ that induces a random walk on $X_n$. One says that there
is a total variation cutoff for the family $(X_n,\pi_n)$ if there
exists a sequence $(t_n)$ of positive reals such that
\begin{enumerate}
\item $\lim_{n \rightarrow \infty} t_n = \infty$;
\item For any $\epsilon \in (0,1)$ and $r_n = [(1+\epsilon)t_n]$, $\lim_{n \rightarrow \infty} ||p_n^{r_n}-\pi_n||=0$;
\item For any $\epsilon \in (0,1)$ and $r_n = [(1-\epsilon)t_n]$,
$\lim_{n \rightarrow \infty} ||p_n^{r_n}-\pi_n||=1$.
\end{enumerate} For the definition of a separation cutoff, one replaces $||p_n^{r_n}-\pi_n||$ by $s(p_n^{r_n},\pi_n)$.

\section{The symmetric group} \label{symirrep}

This section studies the random walk $K$ on $Irr(S_n)$ defined from
the representation $\eta$ whose character is the number of fixed
points. Although not needed for the results in this section, it should
be mentioned that when $Irr(S_n)$ is viewed as the partitions of $n$,
the random walk $K$ has a description in terms of removing and then
reattaching a corner box at each step (see \cite{F6} for a proof).

The primary purpose of this section is to determine the asymptotic
behavior of the separation distance \[ s(r) = \max_{\lambda} \left[ 1 -
\frac{K^r(\hat{1},\lambda)} {\pi(\lambda)} \right].\]

The first step in studying $s(r)$ is determine for which $\lambda$ the
maximum is attained. Part 1 of Lemma \ref{genbound} and Lemma
\ref{diaggroup} imply that \[ \frac{K^r(\hat{1},\lambda)}
{\pi(\lambda)} = \sum_{i=0}^n \left( \frac{i}{n} \right)^r \sum_{g \in
S_n: fp(g) = i} \frac{\chi^{\lambda}(g)}{d_{\lambda}},\] where $fp(g)$
is the number of fixed points of $g$. However since characters can
take both positive and negative values, it is not clear which
$\lambda$ minimizes this expression.

Theorem \ref{symsumpos} will circumvent this difficulty by giving an
expression for $\frac{K^r(\hat{1},\lambda)} {\pi(\lambda)}$ as a sum
of non-negative terms. This result was first derived in our earlier
paper \cite{F5} using the Robinson-Schensted-Knuth correspondence. The
proof presented here is different and uses instead inclusion-exclusion
and the branching rules for the irreducible representations of
$S_n$. As will be seen in Section \ref{glirrep}, it generalizes
perfectly to the group $GL(n,q)$.

As a first step, the following lemma is useful. In its statement, and
throughout this section, we assume familiarity with the concept of
standard tableaux as in Chapters 2 and 3 of \cite{Sag}.

\begin{lemma} \label{sumfix}  Let $d_{\lambda/\mu}$ denote
 the number of standard tableaux of shape $\lambda/\mu$. Then \[ \sum_{g
 \in S_n: fp(g)=i} \chi^{\lambda}(g) = \frac{n!}{i!}
 \sum_{j=0}^{n-i} \frac{(-1)^j}{j!} d_{\lambda/(n-i-j)} .\]
\end{lemma}

\begin{proof} Let $Fix(g)$ denote the set of fixed points of a permutation $g$. Then \begin{eqnarray*} & & \sum_{g \in S_n: fp(g)=i} \chi^{\lambda}(g)\\
 & = & {n \choose i} \sum_{g: Fix(g)=\{n-i+1,\cdots,n\}}
\chi^{\lambda}(g)\\ & = & {n \choose i} \sum_{j=0}^{n-i} (-1)^j
\sum_{A \subseteq \{1,\cdots,n-i\} \atop |A|=j} \sum_{g: Fix(g)
\supseteq A \cup \{n-i+1,\cdots,n\}} \chi^{\lambda}(g)\\ & = & {n
\choose i} \sum_{j=0}^{n-i} (-1)^j {n-i \choose j} \sum_{g: Fix(g)
\supseteq \{n-i-j+1,\cdots,n\}} \chi^{\lambda}(g)\\ & = & {n \choose
i} \sum_{j=0}^{n-i} (-1)^j {n-i \choose j} (n-i-j)! \langle
Res^{S_n}_{S_{n-i-j}}[\chi^{\lambda}],\hat{1} \rangle \\ & = &
\frac{n!}{i!} \sum_{j=0}^{n-i} \frac{(-1)^j}{j!}
d_{\lambda/(n-i-j)}. \end{eqnarray*} The first and third equalities
are since character values are invariant under conjugacy. The second
equality is the inclusion-exclusion principle (Chapter 10 of
\cite{VW}). In the fourth equality,
$Res^{S_n}_{S_{n-i-j}}[\chi^{\lambda}]$ denotes the restriction of
$\chi^{\lambda}$ from $S_n$ to $S_{n-i-j}$. The final equality follows
from the branching rules for irreducible representations of symmetric
groups \cite{Sag}. Note also that when $i=0$, the set
$\{n-i+1,\cdots,n\}$ should be interpreted as the empty
set. \end{proof}

In what follows $P(a,r,n)$ will denote the probability that when $r$
balls are dropped at random into $n$ boxes, there are exactly $a$
occupied boxes. Lemma \ref{combin1} gives an explicit expression for
$P(a,r,n)$. This expression is an exercise on page 103 of \cite{Fe},
but since the proof is simple and motivates an analogous result in
Section \ref{glirrep}, we include it.

\begin{lemma} \label{combin1} (\cite{Fe}) \[ P(a,r,n) =  {n \choose a} \sum_{b=n-a}^n (-1)^{b-(n-a)} {a \choose n-b} \left(1 - \frac{b}{n} \right)^r.\] \end{lemma}

\begin{proof} Clearly $P(a,r,n)$ is ${n \choose a}$ multiplied by the probability that the
 occupied boxes are the first $a$ boxes. By the principle of inclusion
 and exclusion, this is \[ {n \choose a}
 \sum_{s=0}^a (-1)^{a-s} {a \choose s} P_{\leq}(s) \] where
 $P_{\leq}(s)$ is the probability that the set of occupied boxes is
 contained in $\{1,\cdots,s\}$. Noting that $P_{\leq}(s)= \left(
 \frac{s}{n} \right)^r$, the result follows from the change of
 variables $b=n-s$. \end{proof}

Theorem \ref{symsumpos} gives the needed expression for
$\frac{K^r(\hat{1},\lambda)} {\pi(\lambda)}$ as a sum of non-negative
quantities.

\begin{theorem} \label{symsumpos} Let
 $d_{\lambda/\mu}$ denote the number of standard tableaux of shape
 $\lambda/\mu$. Then \[ \frac{K^r(\hat{1},\lambda)} {\pi(\lambda)} =
 \sum_{a=0}^n P(a,r,n)(n-a)! \frac{d_{\lambda/(n-a)}}{d_{\lambda}}.\]
\end{theorem}

\begin{proof} As noted earlier, part 1 of Lemma \ref{genbound} and Lemma
\ref{diaggroup} imply that \[ \frac{K^r(\hat{1},\lambda)}
{\pi(\lambda)} = \sum_{i=0}^n \left( \frac{i}{n} \right)^r \sum_{g \in
S_n: fp(g) = i} \frac{\chi^{\lambda}(g)}{d_{\lambda}}.\] By Lemma
\ref{sumfix} this is \[ \sum_{i=0}^n  \left( \frac{i}{n}
\right)^r \frac{n!}{i!} \sum_{j=0}^{n-i} \frac{(-1)^j}{j!}
\frac{d_{\lambda/(n-i-j)}}{d_{\lambda}}.\] Letting $a=i+j$, this
becomes \begin{eqnarray*} & & \sum_{i=0}^n \sum_{a=i}^n (-1)^{a-i}
\frac{n!}{i!} \frac{1}{(a-i)!} \left( \frac{i}{n} \right)^r
\frac{d_{\lambda/(n-a)}}{d_{\lambda}}\\ & = & \sum_{a=0}^n
\sum_{i=0}^a (-1)^{a-i} \frac{n!}{i!} \frac{1}{(a-i)!}  \left(
\frac{i}{n} \right)^r
\frac{d_{\lambda/(n-a)}}{d_{\lambda}}. \end{eqnarray*} Letting
$b=n-i$, this becomes \[ \sum_{a=0}^n {n \choose a} \sum_{b=n-a}^n (-1)^{b-(n-a)} {a \choose n-b} \left(1 - \frac{b}{n} \right)^r (n-a)!
\frac{d_{\lambda/(n-a)}}{d_{\lambda}}.\] The result follows from Lemma
\ref{combin1}. \end{proof}

\begin{cor} \label{whichmax} The
 quantity $\frac{K^r(\hat{1},\lambda)} {\pi(\lambda)}$ is minimized
 for $\lambda=(1^n)$, corresponding to the sign representation.
\end{cor}

\begin{proof} By Theorem \ref{symsumpos} one wants to find $\lambda$
 minimizing \[ \frac{K^r(\hat{1},\lambda)} {\pi(\lambda)} =
 \sum_{a=0}^n P(a,r,n)(n-a)! \frac{d_{\lambda/(n-a)}}{d_{\lambda}}.\]
 Note that the $a=n-1$ and $a=n$ terms in this expression are
 independent of $\lambda$. Moreover, all other terms are non-negative,
 and vanish when $\lambda$ is the sign representation (which
 corresponds to the partition all of whose parts have size 1). The
 result follows. \end{proof}

Next we use Theorem \ref{symsumpos} and Corollary \ref{whichmax} to
derive both a formula and a precise asymptotic expression for the
separation distance of the Markov chain $K$.

\begin{theorem} \label{asym} Let $s(r)$ be the separation distance between $K^r$ started at the trivial representation and the Plancherel measure $\pi$.
\begin{enumerate}
\item \[ s(r) = \sum_{i=0}^{n-2} (-1)^{n-i} {n \choose i} (n-i-1)
\left( \frac{i}{n} \right)^r .\]
\item For $c$ fixed in $\mathbb{R}$ and $n \rightarrow \infty$, \[ s(n
\log(n)+cn) = 1-e^{-e^{-c}}(1+e^{-c}) + O \left( \frac{\log(n)}{n}
\right).\]
\end{enumerate} \end{theorem}

\begin{proof} (First proof) Theorem \ref{symsumpos} and Corollary \ref{whichmax} imply that  
\begin{eqnarray*} s(r) & = & 1 - \sum_{a=0}^n P(a,r,n) (n-a)! \frac{ d_{(1^n)/(n-a)}}{d_{(1^n)}}\\
& = & 1 - P(n,r,n) - P(n-1,r,n). \end{eqnarray*} By Lemma
\ref{combin1} this is equal to \begin{eqnarray*} & & 1 - \sum_{b=0}^n
(-1)^b {n \choose b} \left(1-\frac{b}{n} \right)^r - n \sum_{b=1}^n
(-1)^{b-1} {n-1 \choose n-b} \left(1-\frac{b}{n} \right)^r\\ & = &
\sum_{b=1}^n (b-1) {n \choose b} (-1)^b \left( 1 - \frac{b}{n}
\right)^r. \end{eqnarray*} Letting $i=n-b$ proves the first assertion.

For the second assertion, we use asymptotics of the coupon collector's
problem: it follows from Section 6 of \cite{CDM} that when $n
\log(n)+cn$ balls are dropped into $n$ boxes, the number of unoccupied
boxes converges to a Poisson distribution with mean $e^{-c}$, and that
the error term in total variation distance is
$O(\frac{\log(n)}{n})$. The chance that a Poisson random variable with
mean $e^{-c}$ takes value not equal to 0 or 1 is
$1-e^{-e^{-c}}(1+e^{-c})$, which completes the proof. \end{proof}

There is a second proof of Theorem \ref{asym}, which uses a connection
with the top to random shuffle of the symmetric group. We prefer the
first proof as the ideas are more elementary (one doesn't need the RSK
correspondence) and generalize to $GL(n,q)$ (see Section
\ref{glirrep}).

\begin{proof} (Second proof) By Corollary \ref{whichmax}, $s(r) = 1 - n!
 K^r(\hat{1},(1^n))$. Theorem 3.1 of \cite{F3} gives that for any
 shape $\lambda$, one has that $K^r(\hat{1},\lambda)$ is equal to the
 chance of obtaining a permutation with Robinson-Schensted-Knuth (RSK)
 shape $\lambda$ after $r$ top to random shuffles started from the
 identity. The only permutation with RSK shape $(1^n)$ is the
 ``longest'' permutation $\pi$, defined by $\pi(i)=n-i+1$ for all
 $i$. It follows from Corollary 2.1 of \cite{DFP} that the separation
 distance for the top to random shuffle is attained at this
 $\pi$. Thus the chain $K$ and the top to random shuffle have the same
 separation distance $s(r)$, so the result follows from page 142 of
 \cite{DFP}. \end{proof}

From the second proof the reader might think that the theory of the
chain $K$ can be entirely understood in terms of the top to random
shuffle. This is untrue. For example if one measures convergence to
the stationary distribution using total variation distance, the top to
random shuffle takes $n \log(n)+cn$ steps to be close to random
\cite{AD1}, but the chain $K$ requires only $\frac{1}{2} n \log(n)+cn$
steps \cite{F6}.

As a final result, we use Corollary \ref{whichmax} and the cycle index
of the symmetric group to give a third proof of part 1 of Theorem
\ref{asym}.

\begin{proof} By Corollary \ref{whichmax}, $s(r) =  1 - \frac{K^r(\hat{1},(1^n))}{\pi(1^n)}$. Part 1 of Lemma \ref{genbound} and Lemma
\ref{diaggroup} imply that \[ \frac{K^r(\hat{1},\lambda)}
{\pi(\lambda)} = \sum_{i=0}^n \left( \frac{i}{n} \right)^r \sum_{g \in
S_n: fp(g) = i} \frac{\chi^{\lambda}(g)}{d_{\lambda}}.\] Specializing
to $\lambda=(1^n)$ implies that \[ s(r) = - \sum_{i=0}^{n-1} \left(
\frac{i}{n} \right)^r \sum_{g \in S_n:fp(g)=i} sign(g).\] Here
$sign(g)=(-1)^{n-c(g)}$, where $c(g)$ is the number of cycles of $g$.

A classic result in combinatorics (see \cite{W}) is the ``cycle
index'' of the symmetric group, which states that \[ 1 + \sum_{n \geq
1} \frac{u^n}{n!} \sum_{g \in S_n} \prod_{j \geq 1} x_j^{n_j(g)} =
\exp \left( \sum_{m \geq 1} \frac{x_m u^m}{m} \right).\] Here $n_j(g)$
is the number of cycles of length $j$ of $g$. Making the substitutions
$x_1=-x$, $x_i=-1$ for $i \geq 2$ and replacing $u$ by $-u$, the cycle
index implies that \begin{eqnarray*} 1 + \sum_{n \geq 1}
\frac{u^n}{n!} \sum_{g \in S_n} sign(g) \cdot x^{fp(g)} & = & \exp \left( xu
- \sum_{m \geq 2} \frac{(-u)^m}{m} \right)\\ & = & \exp(xu-u)
\exp(\log(1+u))\\ & = & \frac{e^{xu}(1+u)}{e^u}. \end{eqnarray*}
Taking the coefficient of $\frac{u^nx^i}{n!}$ on both sides shows that
if $0 \leq i \leq n-1$, then \begin{eqnarray*} \sum_{g \in S_n:
fp(g)=i} sign(g) & = & \frac{n!}{i!}  \left[ \frac{(-1)^{n-i}}{(n-i)!}
+ \frac{(-1)^{n-i-1}}{(n-i-1)!}  \right]\\ & = & (-1)^{n-i+1} {n
\choose i} (n-i-1). \end{eqnarray*} The result now follows from the
previous paragraph. \end{proof}

\section{The general linear group} \label{glirrep}

This section studies random walk on $Irr(GL(n,q))$ in the case that
$\eta$ is the representation of $GL(n,q)$ whose character is
$q^{d(g)}$, where $d(g)$ is the dimension of the fixed space of
$g$. As in Section \ref{symirrep}, we aim to determine the asymptotic
behavior of the separation distance \[ s(r) = \max_{\Lambda} \left[ 1 -
\frac{K^r(\hat{1},\Lambda)}{\pi(\Lambda)} \right].\]

The first step is to find the irreducible representation $\Lambda$ for
which the maximum is attained. Part 1 of Lemma \ref{genbound} and
Lemma \ref{diaggroup} imply that \[
\frac{K^r(\hat{1},\Lambda)}{\pi(\Lambda)} = \sum_{i=0}^n q^{-r(n-i)}
\sum_{g \in GL(n,q):d(g)=i} \frac{\chi^{\Lambda}(g)}{d_{\Lambda}}.
\] Since characters can take both positive and negative values, it is not at all clear
which $\Lambda$ minimizes this expression. As in the symmetric group
case, the key is to find a way to write $\frac{K^r(\hat{1},\Lambda)}
{\pi(\Lambda)}$ as a sum of non-negative terms.

To begin we recall some facts about the representation theory of
$GL(n,q)$. A full treatment of the subject with proofs appears in
\cite{M}, \cite{Z}. As usual a partition
$\lambda=(\lambda_1,\cdots,\lambda_m)$ is identified with its
geometric image $\{(i,j): 1 \leq i \leq m, 1 \leq j \leq \lambda_i \}$
and $|\lambda|=\lambda_1+\cdots+\lambda_m$ is the total number of
boxes. Let $\yy$ denote the set of all partitions, including the empty
partition of size 0.

Given an integer $1 \leq k < n$ and two characters $\chi_1,\chi_2$ of
the groups $GL(k,q)$ and $GL(n-k,q)$, their parabolic induction
$\chi_1 \circ \chi_2$ is the character of $GL(n,q)$ induced from the
parabolic subgroup of elements of the form \[ P = \left \{ \left( \begin{array} {cc} g_1 & * \\
0 & g_2 \end{array} \right) : g_1 \in GL(k,q), g_2 \in GL(n-k,q) \right \} \] by the function $\chi_1(g_1) \chi_2(g_2)$.

A character is called cuspidal if it is not a component of any
parabolic induction. Let ${\it C}_m$ denote the set of cuspidal
characters of $GL(m,q)$ and let ${\it C} = \bigcup_{m \geq 1} {\it C}_m$;
it is known that $|C_m|=\frac{1}{m} \sum_{d|m} \mu(d) (q^{m/d}-1)$ where
$\mu$ is the Moebius function. The unit character of $GL(1,q)$ plays
an important role and will be denoted $e$; it is one of the $q-1$
elements of ${\it C}_1$. Given a family ${\Lambda}: {\it C} \mapsto
\yy$ with finitely many non-empty partitions ${\Lambda}(c)$, its
degree $||\Lambda||$ is defined as $\sum_{m \geq 1} \sum_{c \in {\it
C}_m} m \cdot |{\it \Lambda}(c)|$. A fundamental result is that the
irreducible representations of $GL(n,q)$ are in bijection with the
families of partitions of degree $n$, so we also let $\Lambda$ denote
the corresponding representation.

Let $\vec{e_1},\cdots,\vec{e_n}$ be the standard basis of the vector
space $V$ on which $GL(n,q)$ acts (so the kth component of $\vec{e_j}$
is $\delta_{j,k}$). Define $H(k,q)$ as the subgroup of $GL(n,q)$
consisting of $g$ which fix all of
$\vec{e_1},\cdots,\vec{e_k}$. Equivalently, the elements of $H(k,q)$
are block matrices of the form \[ \left(
\begin{array}{c c} I_k & X \\ 0_{n-k} & Y \end{array} \right) \]
where $I_k$ is a $k$ by $k$ identity matrix, $X$ is any $k$ by $n-k$
matrix with entries in $\mathbb{F}_q$, and $Y$ is any element of
$GL(n-k,q)$. Thus $|H(k,q)|= q^{k(n-k)} |GL(n-k,q)|$. For $\Lambda$ an
element of $Irr(GL(n,q))$, it will be helpful to let
\[ c_k(\Lambda) = \sum_{g \in H(k,q)} \chi^{\Lambda}(g).\] Then $c_k(\Lambda)$ is non-negative, since it is the product of $|H(k,q)|$ and the multiplicity of the trivial representation of $H(k,q)$ in the restricted representation $Res^{GL(n,q)}_{H(k,q)}[\Lambda]$.

It will also be convenient to let $\qb{n}{k}$ denote the q-binomial
coefficient $\frac{(q^n-1) \cdots (q-1)}{(q^k-1) \cdots (q-1)
(q^{n-k}-1) \cdots (q-1)}$, which is equal to the number of $k$
dimensional subspaces of an $n$ dimensional vector space over a finite
field $\mathbb{F}_q$.

\begin{lemma} \label{sumfixq} Let $\Lambda$ be an irreducible representation of $GL(n,q)$. Then \[ \sum_{g \in GL(n,q) \atop d(g)=i} \chi^{\Lambda}(g) = \qb{n}{i} \sum_{j=0}^{n-i} \qb{n-i}{j} (-1)^j q^{{j \choose 2}} c_{i+j}(\Lambda). \]  
\end{lemma}

\begin{proof} Let $Fix(g)$ denote the fixed space of $g$. Also if $A$ is a set of vectors, $\langle A \rangle$ will
 denote their span. Then \begin{eqnarray*} \sum_{g \in GL(n,q) \atop
 d(g)=i} \chi^{\Lambda}(g) & = & \qb{n}{i} \sum_{g \in GL(n,q) \atop
 Fix(g) = \langle \vec{e_1},\cdots,\vec{e_i} \rangle} \chi^{\Lambda}(g)\\ & = &
 \qb{n}{i} \sum_{j=0}^{n-i} \sum_{W \supseteq
 \langle \vec{e_1},\cdots,\vec{e_i} \rangle \atop dim(W)=i+j} (-1)^j q^{{j \choose
 2}} \sum_{g \in GL(n,q) \atop Fix(g) \supseteq W} \chi^{\Lambda}(g)\\
 & = & \qb{n}{i} \sum_{j=0}^{n-i} \qb{n-i}{j} (-1)^j q^{{j \choose 2}}
 \sum_{g \in GL(n,q) \atop Fix(g) \supseteq
 \langle \vec{e_1},\cdots,\vec{e_{i+j}} \rangle} \chi^{\Lambda}(g)\\ & = & \qb{n}{i}
 \sum_{j=0}^{n-i} \qb{n-i}{j} (-1)^j q^{{j \choose 2}}
 c_{i+j}(\Lambda). \end{eqnarray*}

  The first and third equalities used the fact that if $W_1,W_2$ are
subspaces of $V$ of equal dimension, then $\{ g:Fix(g)=W_1 \}$ and
$\{g:Fix(g)=W_2 \}$ are conjugate in $GL(n,q)$. The second equality
used Moebius inversion on the lattice of subspaces of a vector space
(Chapter 25 of the text \cite{VW}). The third equality also used the
fact that the number of $i+j$ dimensional subspaces of $V$ containing
$\langle \vec{e_1},\cdots,\vec{e_i} \rangle$ is
$\qb{n-i}{j}$. \end{proof}

In what follows we let $P_q(a,r,n)$ denote the probability that the
span of $\vec{v_1},\cdots,\vec{v_r}$ is $a$ dimensional, where the
$r$ vectors are chosen uniformly at random from an n-dimensional
vector space over $\mathbb{F}_q$. Lemma \ref{countsub} gives a formula
for $P_q(a,r,n)$.

\begin{lemma} \label{countsub} \[ P_q(a,r,n) = \qb{n}{a} \sum_{b=n-a}^{n} (-1)^{b-(n-a)} q^{{b-(n-a) \choose 2}} \qb{a}{n-b} q^{-rb}.\] \end{lemma}

\begin{proof} Clearly $P_q(a,r,n)$ is $\qb{n}{a}$ multiplied by the
 chance that the span of $\vec{v_1},\cdots,\vec{v_r}$ is exactly the
 $a$ dimensional subspace consisting of vectors whose last $n-a$
 coordinates are 0. One applies Moebius inversion on the lattice of
 subspaces of an n-dimensional vector space over $\mathbb{F}_q$
 (Chapter 25 of \cite{VW}) to conclude that \[ P_q(a,r,n) = \qb{n}{a}
 \sum_{s=0}^a \sum_{W \subseteq  \langle \vec{e_1},\cdots,\vec{e_a} \rangle \atop
 dim(W)=s} (-1)^{a-s} q^{{a-s \choose 2}} P_{\leq(W)}. \] Here
 $\vec{e_1},\cdots,\vec{e_n}$ is the standard basis of $V$ and
 $P_{\leq}(W)$ is the probability that the span of
 $\vec{v_1},\cdots,\vec{v_r}$ is contained in $W$. Clearly
 $P_{\leq}(W) = q^{-r(n-dim(W))}$. Thus \[ P_q(a,r,n) = \qb{n}{a}
 \sum_{s=0}^a \qb{a}{s} (-1)^{a-s} q^{{a-s \choose 2}} q^{-r(n-s)}, \]
 and the result follows by the change of variables
 $b=n-s$. \end{proof}

Theorem \ref{sympsumposq} is a key result of this section; it
expresses $\frac{K^r(\hat{1},\Lambda)}{\pi(\Lambda)}$ as a sum of
non-negative terms.

\begin{theorem} \label{sympsumposq} \[ \frac{K^r(\hat{1},\Lambda)}{\pi(\Lambda)} =
 \sum_{a=0}^n P_q(a,r,n) \frac{c_a(\Lambda)}{d_{\Lambda}}.\]
\end{theorem}

\begin{proof} Part 1 of Lemma \ref{genbound} and Lemma \ref{diaggroup} imply that
 \[ \frac{K^r(\hat{1},\Lambda)}{\pi(\Lambda)} = \sum_{i=0}^n
 q^{-r(n-i)} \sum_{g \in GL(n,q):d(g)=i}
 \frac{\chi^{\Lambda}(g)}{d_{\Lambda}}.
\] By Lemma \ref{sumfixq}, this is \[ \sum_{i=0}^n
 q^{-r(n-i)} \qb{n}{i} \sum_{j=0}^{n-i} \qb{n-i}{j} (-1)^j q^{{j
 \choose 2}} \frac{c_{i+j}(\Lambda)}{d_{\Lambda}}.\] Letting $a=i+j$,
 this becomes
 \begin{eqnarray*} & & \sum_{i=0}^n q^{-r(n-i)} \qb{n}{i} \sum_{a=i}^n
 \qb{n-i}{a-i}(-1)^{a-i} q^{{a-i \choose 2}}
 \frac{c_{a}(\Lambda)}{d_{\Lambda}} \\ & = & \sum_{a=0}^n
 \frac{c_{a}(\Lambda)}{d_{\Lambda}} \sum_{i=0}^a q^{-r(n-i)} \qb{n}{i}
 \qb{n-i}{a-i}(-1)^{a-i} q^{{a-i \choose 2}}. \end{eqnarray*} Setting
 $b=n-i$ this becomes \begin{eqnarray*} & & \sum_{a=0}^n
 \frac{c_{a}(\Lambda)}{d_{\Lambda}} \sum_{b=n-a}^n q^{-rb} \qb{n}{b}
 \qb{b}{a-(n-b)}(-1)^{a-(n-b)} q^{{a-(n-b) \choose 2}}\\ & = &
 \sum_{a=0}^n \frac{c_a(\Lambda)}{d_{\Lambda}} \qb{n}{a}
 \sum_{b=n-a}^{n} (-1)^{b-(n-a)} q^{{b-(n-a) \choose 2}} \qb{a}{n-b}
 q^{-rb}. \end{eqnarray*} The result now follows from Lemma
 \ref{countsub}. \end{proof}

\begin{cor}
 \label{whichmaxq} Suppose that $GL(n,q) \neq GL(1,2)$.
\begin{enumerate}
\item The quantity $\frac{K^r(\hat{1},\Lambda)}{\pi(\Lambda)}$ is
minimized for any $\Lambda$ which satisfies $\Lambda(e)=\emptyset$,
and such $\Lambda$ exist.
\item The separation distance $s(r)$ of $K$ started at the trivial
representation is \[ 1 - P_q(n,r,n) = \sum_{b=1}^n (-1)^{b+1} q^{{b \choose 2}} \qb{n}{b} q^{-rb}.\]
\end{enumerate}
\end{cor}

\begin{proof} The formula for $|C_m|$ stated earlier in this section implies the existence of $\Lambda$ with $\Lambda(e)=\emptyset$. Proposition 5.4 of \cite{F6} states that for any $\Lambda$, if
 $K^r(\hat{1},\Lambda)>0$, then the largest part of $\Lambda(e)$ is at
 least $n-r$. It follows that if $\Lambda(e)=\emptyset$ then
 $K^{n-1}(\hat{1},\Lambda)=0$. By Theorem \ref{sympsumposq} and the
 fact that $P_q(a,n-1,n)>0$ for $0 \leq a \leq n-1$, it follows that
 if $\Lambda(e)=\emptyset$ then $c_a(\Lambda)=0$ for $0 \leq a \leq
 n-1$. Thus if $\Lambda(e)=\emptyset$, only the $a=n$ term in Theorem
 \ref{sympsumposq} can be non-vanishing, but this term is independent
 of $\Lambda$, since $\frac{c_n(\Lambda)}{d_{\Lambda}}=1$ for all
 $\Lambda$. This implies the first part of the lemma since all terms
 in Theorem \ref{sympsumposq} are nonnegative. It also implies that
 $s(r)=1-P_q(n,r,n)$, and the equality in the second assertion follows
 from Lemma \ref{countsub}. \end{proof}

Theorem \ref{asympq} bounds the separation distance $s(r)$ and
determines its exact asymptotic behavior.

\begin{theorem} \label{asympq} Suppose that $GL(n,q) \neq GL(1,2)$.
\begin{enumerate}
\item If $r<n$, then $s(r)=1$. If $c \geq 0$, then \[ \frac{1}{q^{c+1}}-\frac{4}{q^{2c+3}} \leq s(n+c) \leq \frac{2}{q^{c+1}}.\]
\item Let $c \geq 0$ be fixed. Then \[ \lim_{n \rightarrow \infty} s(n+c) = 1 -
\prod_{m=1}^{\infty} ( 1 - q^{-(c+m)}).\]
\end{enumerate}
\end{theorem}

\begin{proof} The first two sentences of the proof of Corollary \ref{whichmaxq} implies that if $r<n$, then $s(r)=1$. To upper bound $s(n+c)$, one checks that since $q \geq 2$, the sum in the second part of Corollary \ref{whichmaxq} is alternating with terms of decreasing magnitude. Thus the sum is upper bounded by its first term, $\qb{n}{1} q^{-(n+c)}$, which is easily seen to be less than $2q^{-(c+1)}$, since $q \geq 2$. For the lower bound, note that the first term in the second part of Corollary \ref{whichmaxq} is at least $q^{-(c+1)}$ and that the second term $-\frac{\qb{n}{2}}{q^{2n+2c-1}}$ is at least $-4q^{-(2c+3)}$ since $q \geq 2$. This proves the first part of the theorem.

To prove the second part, rewrite the expression for $s(n+c)$ in part
2 of Corollary \ref{whichmaxq} as \[ - \sum_{b=1}^n \frac{(-1)^b
(1-1/q^n)(1/q-1/q^n) \cdots (1/q^{b-1}-1/q^n)}{q^{(c+1)b} (1-1/q)
\cdots (1-1/q^b)}.\] It is straightforward to see that as $n
\rightarrow \infty$ this converges to \[ - \sum_{b=1}^{\infty}
\frac{(-1)^b}{q^{(c+1)b} q^{{b \choose 2}} (1-1/q) \cdots
(1-1/q^b)}.\] An identity of Euler (Corollary 2.2 in \cite{An}) states
that \[ 1 + \sum_{b=1}^{\infty} \frac{t^b}{q^{{b \choose 2}} (1-1/q)
\cdots (1-1/q^b)} = \prod_{m=0}^{\infty} (1+tq^{-m}) \] for
$|t|<1,|q|>1$. The result follows by applying this identity with
$t=-1/q^{(c+1)}$. \end{proof}

\section{Lagrange-Sylvester Interpolation} \label{lagrange}

For certain one dimensional random walks, namely stochastically
monotone birth-death chains, Lagrange-Sylvester interpolation allows
one to study separation distance knowing only the eigenvalues (and not
the eigenvectors) of the Markov chain; see \cite{Br}, \cite{DSa} and
the remarks following Proposition \ref{noeigenvec}. The purpose of
this section is to give examples of higher dimensional state spaces
(namely $Irr(S_n)$ and $Irr(GL(n,q))$) where the methodology is
useful.

To begin we review the Lagrange-Sylvester interpolation approach to
diagonalizable matrices, and hence to reversible Markov chains. A
textbook discussion in the matrix setting appears in \cite{Ga}, and
the paper \cite{Br} uses the language of Markov chains.

As usual, $K$ is a Markov chain on a finite set $X$ and is
reversible with respect to a distribution $\pi$. If $\pi(x)>0$ for all
$x$, then letting $A$ be a diagonal matrix whose $(x,x)$ entry is
$\pi(x)$, it follows that $A^{1/2} K A^{-1/2}$ is symmetric. Hence $K$
is conjugate to a diagonal matrix $D$, whose entries $d_1,\cdots,d_n$
are the eigenvalues of $K$. Thus if $f,g$ are polynomials with
$f(d_i)=g(d_i)$ for $i=1,\cdots,n$ then $f(K)=g(K)$.

Let $\lambda_1,\cdots,\lambda_m$ be the distinct eigenvalues of $K$
(so $m \leq |X|$). Define $g_r(s) = s^r$ and \[ f_r(s) = \sum_{i=1}^m
\lambda_i^r \left[ \prod_{j \neq i}
\frac{s-\lambda_j}{\lambda_i-\lambda_j} \right] .\] Since
$f_r(\lambda_i)=g_r(\lambda_i)$ for $i=1,\cdots,m$, it follows that
$f_r(D)=g_r(D)$. Thus $f_r(K)=g_r(K)$ which gives that \[ K^r =
\sum_{i=1}^m \lambda_i^r \left[ \prod_{j \neq i}
\frac{K-\lambda_j I}{\lambda_i-\lambda_j} \right] .\] As noted in
\cite{Br}, expanding this expresses $K^r$ in terms of
$I,K,\cdots,K^{m-1}$ as follows: \[ K^r = \sum_{a=1}^m K^{a-1}
(-1)^{m-a} \sum_{i=1}^m \lambda_i^r \prod_{j \neq i} (\lambda_i
- \lambda_j)^{-1} \sum_{\alpha \in c(m,i,m-a)} ( \prod_{s \in
\alpha} \lambda_s ). \] Here $c(m,i,m-a)$ consists of the ${m-1
\choose a-1}$ subsets of size $m-a$ from $\{j: 1 \leq j \leq m, j \neq
i \}$.

For the next proposition it is useful to define the distance
$dist(x,y)$ between $x,y \in X$ as the smallest $r$ such that
$K^r(x,y)>0$. For the special case of birth-death chains on the set
$\{0,1,\cdots,d\}$, Proposition \ref{noeigenvec} appears in
\cite{DF} and \cite{Br}.

\begin{prop} \label{noeigenvec}
 Let $K$ be a reversible ergodic Markov on a finite set $X$. Let
 $1,\lambda_1,\cdots,\lambda_d$ be the distinct eigenvalues of $K$ (so
 $d+1 \leq |X|$). Suppose that $x,y$ are elements of $X$ with
 $dist(x,y)=d$. Then for all $r \geq 0$,
 \[ 1 - \frac{K^r(x,y)}{\pi(y)} =  \sum_{i=1}^d \lambda_i^r
\left[ \prod_{j \neq i}
\frac{1-\lambda_j}{\lambda_i-\lambda_j} \right].\]
 \end{prop}

\begin{proof} Since $dist(x,y)=d$, the Lagrange-Sylvester expansion of $K^r$ in terms of $I,K,\cdots,K^d$ gives that
 \[ K^r(x,y) = K^d(x,y) \left(  \prod_{j} (1 - \lambda_j)^{-1} - \sum_{i=1}^d  \lambda_i^r (1-\lambda_i)^{-1} \prod_{j \neq i} (\lambda_i - \lambda_j)^{-1} \right).\] By Lemma \ref{genbound}, ergodicity of $K$ implies that $\pi(y) =
lim_{r \rightarrow \infty} K^r(x,y)$. Thus \[ \pi(y) = K^{d}(x,y) \prod_{j}
(1-\lambda_j)^{-1}, \] which implies that \[ \frac{K^r(x,y)}{\pi(y)} = 1 -
\sum_{i=1}^d \lambda_i^r \left[ \prod_{j \neq i}
\frac{1-\lambda_j}{\lambda_i-\lambda_j} \right]. \] \end{proof}

{\bf Remarks:}

\begin{enumerate}
\item Let $K$ be a reversible Markov chain on a finite set $X$. As in
 Proposition \ref{noeigenvec}, let $1,\lambda_1,\cdots,\lambda_d$ be
 the distinct eigenvalues of $K$ (so $d+1 \leq |X|$). From the
 expansion of $K^r$ in terms of $I,K,\cdots,K^d$ it follows that
 $dist(x,y) \leq d$ for any states of the chain. Thus the hypothesis
 of Proposition \ref{noeigenvec} is an extremal case.

\item Let $K$ be a birth-death chain on the set $\{0,\cdots,d\}$, with
transition probabilities \begin{eqnarray*} a_x & = & K(x,x-1) , \
x=1,\cdots, d \\ b_x & = & K(x,x), \ x=0, \cdots, d \\ c_x & = &
K(x,x+1), \ x=0, \cdots, d-1. \end{eqnarray*} Suppose that $a_x>0$ for
$0<x \leq d$ and that $c_x>0$ for $0 \leq x<d$. Such chains are
reversible with respect to the stationary distribution \[ \pi(x) = Z
\prod_{i=1}^x \frac{c_{i-1}}{a_i}, \] where $Z$ is a normalizing
constant. Supposing further that $c_x+a_{x+1} \leq 1$ for $0 \leq x<
d$ (such chains are called monotone chains), then \cite{DF} showed
that the separation distance for the chain started at 0 is equal to \[
s(r) = 1 - \frac{K^r(0,d)}{\pi(d)}.\] Applying Proposition
\ref{noeigenvec} with $x=0,y=d$ one recovers the lovely result of
Diaconis and Fill \cite{DF} expressing the separation distance
entirely in terms of the eigenvalues of $K$: \[ s(r) = \sum_{i=1}^d
\lambda_i^r \left[ \prod_{j \neq i} \frac{1-\lambda_j}
{\lambda_i-\lambda_j} \right]. \] This fact was used in \cite{DSa} to
give a necessary and sufficient spectral condition for the existence
of a separation cutoff for monotone birth death chains.

\end{enumerate}

Next we use Proposition \ref{noeigenvec} to give eigenvector-free
proofs of the formulas for separation distance for the random walks on
$Irr(S_n)$ and $Irr(GL(n,q))$ analyzed in Sections \ref{symirrep} and
\ref{glirrep}. It should be emphasized that as in the proofs of
Sections \ref{symirrep} and \ref{glirrep}, one still needs to know at
what representation the separation distance is attained.

\begin{proof} (Fourth proof of part 1 of Theorem \ref{asym}) By
 Corollary \ref{whichmax}, the separation distance is equal to \[ s(r)
 = 1 - \frac{K^r(\hat{1},(1^n))}{\pi((1^n))}, \] where $(1^n)$ is the
 partition corresponding to the sign representation. As was
 mentioned at the beginning of Section \ref{symirrep} and proved in
 \cite{F6}, the Markov chain $K$ has a description as a random walk on
 partitions in which one removes and adds a box at each step. From
 that description it is clear that the trivial representation
 (corresponding to the partition $(n)$) and the sign representation
 $(1^n)$ are distance $n-1$ apart. By part 1 of Lemma \ref{diaggroup},
 the chain $K$ has $n$ distinct eigenvalues, namely $\frac{i}{n}$
 where $0 \leq i \leq n-2$ or $i=n$. Thus Proposition \ref{noeigenvec}
 implies that \begin{eqnarray*} s(r) & = & \sum_{i=0}^{n-2} \left(
 \frac{i}{n} \right)^r \prod_{j \neq i \atop 0 \leq j \leq n-2}
 \frac{\left( 1-\frac{j}{n} \right)}{\left( \frac{i}{n} - \frac{j}{n}
 \right)} \\
& = &  \sum_{i=0}^{n-2} \left( \frac{i}{n} \right)^r \prod_{j \neq i \atop 0 \leq j \leq n-2} \frac{n-j}{i-j} \\
& = &  \sum_{i=0}^{n-2} \left( \frac{i}{n} \right)^r \frac{n!}{n-i} \prod_{j \neq i \atop 0 \leq j \leq n-2} \frac{1}{i-j} \\
& = & \sum_{i=0}^{n-2} (-1)^{n-i} {n \choose i} (n-i-1)  \left( \frac{i}{n} \right)^r. \end{eqnarray*} \end{proof}

A similar argument works for the general linear case.

\begin{proof} (Second proof of part 2 of Corollary \ref{whichmaxq}) By part 1 of
 Corollary \ref{whichmaxq}, the separation distance is equal to \[
 s(r) = 1 - \frac{K^r(\hat{1},\Lambda)}{\pi(\Lambda)}, \] where
 $\Lambda$ is any representation satisfying $\Lambda(e)=\emptyset$. By
 part 1 of Lemma \ref{diaggroup}, the chain $K$ has $n+1$ distinct
 eigenvalues, namely $q^{-i}$ for $0 \leq i \leq n$. By the first
 remark after Proposition \ref{noeigenvec}, this implies that
 $dist(\hat{1},\Lambda) \leq n$. Proposition 5.4 of \cite{F6} states
 that for any $\Lambda$, if $K^r(\hat{1},\Lambda)>0$ then the largest
 part of $\Lambda(e)$ is at least $n-r$; it follows that
 $dist(\hat{1},\Lambda) \geq n$. Thus $dist(\hat{1},\Lambda)=n$, and
 so Proposition \ref{noeigenvec} can be applied with
 $x=\hat{1},y=\Lambda$. Using the notation $(1/q)_k = (1-1/q) \cdots
 (1-1/q^k)$, one obtains that \begin{eqnarray*} s(r) & = &
 \sum_{i=1}^n q^{-ir} \prod_{j \neq i \atop 1 \leq j \leq n} \frac{1-q^{-j}}{q^{-i}-q^{-j}} \\
 & = & \sum_{i=1}^n q^{-ir} \frac{(1/q)_n}{(1-q^{-i})} \prod_{j=1}^{i-1}
 \frac{1}{q^{-i}-q^{-j}} \prod_{j=i+1}^n \frac{1}{q^{-i}-q^{-j}}\\ & = &
 \sum_{i=1}^n q^{-ir} \frac{(1/q)_n}{(1-q^{-i})} \frac{(-1)^{i-1} q^{{i
 \choose 2}}}{(1/q)_{i-1}} \frac{q^{i(n-i)}}{(1/q)_{n-i}}\\ & = &
 \sum_{i=1}^n (-1)^{i+1} q^{{i \choose 2}} \qb{n}{i}
 q^{-ir}. \end{eqnarray*}
 \end{proof}

\end{document}